\documentclass[12pt]{article} \textwidth=6in \oddsidemargin=0in
\usepackage{amssymb} 
\begin{document}
\def\ov{\over} \def\d{\delta} \def\ld{\ldots}
\def\ee{\end{equation}}\def\({\left(} \def\){\right)}
\def\lb{\Big[} \def\rb{\Big]} \def\iy{\infty} \def\be{\begin{equation}}
\def\ee{\end{equation}} \def\Z{\mathbb Z} \def\cd{\cdots}
\def\P{\mathbb P} \def\n{\nu} \def\s{\sigma} \def\x{\xi}
\def\ch{\raisebox{.3ex}{$\chi$}} \def\p{\pi} \def\a{\alpha}
\def\X{\mathcal X} \def\Y{\mathcal Y} \def\Tp{{T_i\p}} \def\Ts{{T_i\s}}
\newcommand{\xs}[1]{\x_{\s(#1)}} \def\inv{^{-1}} \def\si{\s\inv}
\def\S{\mathbb S} \def\bc{\begin{center}} \def\ec{\end{center}}
\def\a{\alpha} \def\b{\beta} \newcommand{\C}[1]{{\cal C}_{#1}}
\def\ep{\varepsilon} \def\ph{\varphi} \def\noi{\noindent} \def\l{\ell}
\def\e{\eta} \newcommand{\vs}[1]{\vskip#1ex} \def\t{\tau}
\def\A{\mathcal A} \def\H{\mathcal H} \def\sss{\s^{1\,2}}
\def\i{\iota} \def\l{\ell} \def\ep{\varepsilon} \def\c{\circ} \def\noi{\noindent} \def\cdt{\,\cdot\,} \def\F{\mathcal{F}} \def\S{\mathcal{S}} \def\I{\mathcal I} \def\ti{\times} \def\m{\mu} \def\pii{\p\inv} \def\r{\rho} 

\begin{center}{\Large\bf On the Asymmetric Simple Exclusion Process\\ 
\vskip1ex with Multiple Species}\end{center}

\bc{\large\bf Craig A.~Tracy}\\
{\it Department of Mathematics \\
University of California\\
Davis, CA 95616, USA\\
email:\,tracy@math.ucdavis.edu}\ec
 
\bc{\large \bf Harold Widom}\\
{\it Department of Mathematics\\
University of California\\
Santa Cruz, CA 95064, USA\\
email:\,widom@ucsc.edu}\ec  

\bc {\bf I. Introduction}\ec

The one-dimensional asymmetric simple exclusion process (ASEP) \cite{Li1, Li2, GM} is one of the simplest interacting particle systems with a single conservation law (density of particles) and as such is a basic model in both probability theory and nonequilibrium statistical physics.  Its importance is further enhanced by the fact that weakly asymmetric limits of ASEP
distributions can be interpreted as distributions of the height function which solves the Kardar-Parisi-Zhang equation \cite{KPZ, ACQ, SS1}. 

Recall that in the asymmetric simple exclusion process  particles are at sites of the lattice $\Z$.\footnote{Many authors consider ASEP on the circle or the lattice $[1,L]$ with open boundary conditions.} Each particle waits exponential time, then with probability $p$ it moves one step to the right if the site is unoccupied, otherwise it stays put; and with probability $q=1-p$ it moves one step to the left if the site is unoccupied, otherwise it stays put. Each particle does this independently of the other particles.  For a finite number of particles this defines a Markov process; and for infinitely many particles, with further work \cite{Li1} this too defines a Markov process.

In \textit{multispecies ASEP} particles belong to different species, labelled $1, 2, \ld, M$. Particles of a higher species have priority over those of a lower species.\footnote{This is sometimes called the $M+1$ species model, empty sites behaving as particles of another species. With our convention, a particle of species $M$ is first-class, having priority over all others.}  Thus, if a particle of species $s$ tries to move to a neighboring site occupied by a particle of species $s'$ it is blocked if $s\le s'$, but if $s>s'$ the particles interchange positions.  Second-class particles were introduced by Liggett \cite{Li3} and subsequently developed and generalized by several authors
 \cite{AB, AAMP, AKSS,C,DE,FFK,W-K}.
 
A configuration in ASEP with $N$ particles is the set of occupied sites
\[X=\{x_1,\ld,x_N\},\ \ \ (x_1<\cd<x_N).\]
Theorem 2.1 of \cite{tw1} (with proof corrected in \cite{tw2}) was a formula for $P_Y(X;\,t)$, the probability that the system is in configuration $X$ at time $t$, given that the initial configuration was $Y=\{y_1,\ld,y_N\}$. It is a sum over the permuation group $\S_N$ of multiple integrals. If $p\ne 0$ then
\be P_Y(X;\,t)=\sum_{\s\in \S_N} {1\ov(2\pi i)^N}\int_{\C{r}^N} A_\s(\x)\,\prod_i\xs{i}^{x_i}\,\prod_i \Big(\x_i^{-y_i-1}\,e^{\ep(\x_i)\,t})\;d^N\x,\label{prob1}\ee
where $\C{r}$ is a circle about zero in $\mathbb C$ with sufficiently small radius $r$, where
\[\ep(\x)=p\,\x\inv+q\,\x-1,\]
and where $A_\s(\x_1,\ld,\x_N)$ is given explicitly by (\ref{A}) below.

In multispecies ASEP a configuration $\X$ is a pair 
$(X,\,\p)$ where $X=\{x_1,\ld,x_N\}$ as before and $\p$ is a function from $[1,\,N]$ to $[1,\,M]$. If the system is in configuration $\X$ then the $i$th particle from the left is at $x_i$ and belongs to species $\p_i$.
A special case is that of first- and second-class particles, a first-class particle having priority over a second-class particle. For example, if $\p=(1\ 2\ 2\ 2)$ the left-most particle is second-class and the other three are first-class. 

The purpose of this paper is to establish for multispecies ASEP a  formula analogous to (\ref{prob1}) for  $P_\Y(\X;\,t)$, the probability that the system is in configuration $\X=(X,\,\p)$ at time $t$, given that the initial configuration is $\Y=(Y,\,\n)$. We show that there is an entirely analogous formula,
\be P_\Y(\X;\,t)=\sum_{\s\in \S_N} {1\ov(2\pi i)^N}\int_{\C{r}^N} A_\s^\p(\x)\,\prod_i\xs{i}^{x_i}\,\prod_i \Big(\x_i^{-y_i-1}\,e^{\ep(\x_i)\,t}\Big)\;d^N\x,\label{prob2}\ee
but now the factors $A_\s^\p$ are not (except in special cases) given explicitly. They are determined by (\ref{Ah}) and (\ref{h}) below. Note that they also depend on $\n$.

In the next section we present the proof of (\ref{prob1}) in some detail. This is partly because for a correct proof one must refer to both \cite{tw1} and \cite{tw2}, but also because we shall show how to prove (\ref{prob2}) by a modification of the proof of (\ref{prob1}). 

There is another difference here. The factors $A_\s^\p$ must satisfy the family of identities determined by (\ref{Ah}) and (\ref{h}) below. For ordinary ASEP there is no issue about the existence of a solution since it can be written down explicitly. But for multispecies ASEP we must show that the identities define $A_\sigma^\pi$ consistently.  If, as in \cite{AB}, the multispecies ASEP is formulated as a nested Bethe Ansatz problem, one is led to show that the \textit{Yang-Baxter equations} \cite{B,Y} are satisfied.  That the model has nontrivial
solutions to the  Yang-Baxter equations can be traced back to early work of Perk and Schultz on multistate vertex models \cite{PS}.  It was
Alcaraz \textit{et al.} \cite{ADHR} who recognized that these Bethe Ansatz solvable multistate vertex models lead to integrable stochastic models.  
In our formulation below, these consistency conditions are stated in terms of representations of braid relations.

In the last section we focus on the explicit determination of some $A_\s^\p$, with particular attention to the case of one second-class particle.

\bc{\bf II. ASEP with One Species} \ec

\bc{\bf A. Bethe Ansatz solution}\ec

The probability $P_Y(X;\,t)$ satisfies the differential equation (the {\it master equation}) 
\[{\partial u\ov\partial t}=\sum_{i=1}^N
\lb p\,u(x_i-1)\,(1-\d(x_i-x_{i-1}-1))
+q\,u(x_i+1)\,(1-\d(x_{i+1}-x_i-1))\]
\be -p\,u(x_i)\,(1-\d(x_{i+1}-x_i-1))-q\,u(x_i)\,(1-\d(x_i-x_{i-1}-1))\rb.\label{eq1}\ee
Here $u$ is $u(x_1,\ld,x_N;\,t)$, and in the $i$th summand above the $j$th variable is $x_j$ when $j\ne i$. We use the convention that any $\d$ term involving $x_0$ or $x_{N+1}$ is zero.

The probability $P_Y(X;\,t)$ is the solution of this equation that also satisfies the initial condition
\be u(X;\,0)=\d_Y(X).\label{ic}\ee

The particles interact through the exclusion constraint. If they did not interact then $P_Y(X;t)$ would satisfy the differential equation 
\be {\partial u\ov\partial t}=\sum_{i=1}^N
\lb p\,u(x_i-1)+q\,u(x_i+1)-u(x_i)\rb.\label{eq2}\ee
 If $u$ satisfies this equation then it would also satisfy (\ref{eq1}) if in addition the difference of the right sides were zero. This difference is
\[\sum_{i=1}^{N-1}\lb p\,u(x_i,\,x_{i+1}-1)+q\,u(x_i+1,\,x_{i+1})-
u(x_i,\,x_{i+1})\rb\,\d(x_{i+1}-x_i-1).\]
Here we displayed entries $i$ and $i+1$ in $u(x)$. For the first and last summands, we changed indices of summation from what they were in (\ref{eq1}). 

Thus $u(x)$ satisfies (\ref{eq1}) if it satisfies both (\ref{eq2}) and, for $i=1,\ld,N-1$, the {\it boundary conditions} 
\be p\,u(x_i,x_i)+q\,u(x_i+1,x_i+1)-u(x_i,x_i+1)=0.\label{bc}\ee

For any nonzero complex numbers $\x_1,\ld,\x_N$, a solution of (\ref{eq2}) is
$\prod_i\(\x_i^{x_i}\,e^{\ep(\x_i)\,t}\)$.
We may permute the $\x_i$, take linear combinantions, and integrate.  In this way we obtain a class of solutions  
\[\int \sum_{\s\in \S_N} F_\s(\x)\,\prod_i\xs{i}^{x_i}\,\prod_i e^{\ep(\x_i)\,t}\,d^N\x,\]
where the functions $F_\s(\x)$ are arbitrary. We look for $F_\s$ such that the integrand satisfies the boundary conditions pointwise. This is the {\it Bethe Ansatz}.

Substituting the left side of (\ref{bc}) into the part of the integrand that depends on $x_i$ we get
$(\xs{i}\,\xs{i+1})^{x_i}$ times
\[(p+q\,\xs{i}\,\xs{i+1}-\xs{i+1})\,F_\s.\]

Define $\Ts$ to be $\s$ with the entries $\s(i)$ and $\s(i+1)$ interchanged. If we replace $\s$ by $\Ts$ nothing else in the integrand changes, so a sufficient condition that the integrand is zero is that
\[(p+q\,\xs{i}\,\xs{i+1}-\xs{i+1})\,F_\s+
(p+q\,\xs{i}\,\xs{i+1}-\xs{i})\,F_{\Ts}=0.\]
If we define
\[S(\x,\x')=-{p+q\x\x'-\x\ov p+q\x\x'-\x'},\]
then the conditions become
\[{F_\Ts\ov F_\s}=S(\xs{i+1},\xs{i}).\]

We can find the general solution of this system of equations for the $F_\s$. An {\it inversion} in $\s$ is a pair $(i,j)$ with $i>j$ and 
$\si(i)<\si(j)$. It is straightforward to check that one solution is
\be A_\s(\x)=\prod_{{\rm inversions}\,(i,j)}S(\x_i,\x_j).\label{A}\ee
All $F_\s$ are determined by $F_{e}$, where $e$ is the identity permutation. Since this can be an arbitrary function $\ph(\x)$, 
the general solution is $F_\s(\x)=A_\s(\x)\;\ph(\x)$.
\pagebreak

\bc{\bf B. Satisfying the initial condition}\ec
 
We choose $\ph(\x)$ so that the initial condition is satisfied by the $\s=e$ summand. Since
\[{1\ov(2\pi i)^N}\int_{\C{}^N}\, \prod_i\x_i^{x_i-y_i-1}\,d^N\x=\d_Y(X),\]
where $\C{}$ is a circle about zero, we take
$\ph(\x)=(2\pi i)^{-N}\prod_i\x_i^{-y_i-1}$ and the domain of integration to be ${\C{}^N}$. 

When $\s\ne e$ it matters which contours of integration we take  because of the denominators in the $A_\s$. If $p\ne0$ then all the denominators will be nonzero on and inside $\C{r}$ if $r$ is small enough. It is such an $r$ that we take in (\ref{prob1}). 

Denote by $I(\s)$ the $\s$-summand in (\ref{prob1}) with $t=0$. To prove (\ref{prob1}) we must show that
\[\sum_{\s\ne e}I(\s)=0. \]

This will be proved by induction on $N$. For $N=1$ there is nothing to show. Assuming the result for $N-1$, the sum over all permutations in $\S_N\backslash\{e\}$ such that $\s(N)=N$  will be zero, by the induction hypothesis. (No $S$-factor in $A_\s$ will involve $\x_N$, and there is an obvious correspondence between $\s\in\S_{N-1}$ and those $\s\in\S_N$ satisfying $\s(N)=N$.) So it suffices to show that
\be\sum_{\s(N)\ne N}I(\s)=0.\label{Isum}\ee

For each nonempty subset $B$ of $[1,\,N-1]$ define
\[\S_N(B)=\{\s\in\S_N:\textrm{the inversions in $\s$ involving $N$ are the $(N,i)$ with}\ i\in B\}.\]
We shall show that for each $B$ we have
\be\sum_{\s\in\S_N(B)}I(\s)=0.\label{IBsum}\ee
Once we have this, (\ref{Isum}) will follow since the set of $\s$ with $\s(N)\ne N$ is the disjoint union of the various $\S_N(B)$.

\noi{\bf Start of the proof}. The integrands in $I(\s)$ with $\s\in\S_N(B)$ may be written
\[\prod_{i\in B}S(\x_N,\x_i)\times
\prod_{i\le N}\x_i^{x_{\s\inv(i)}-y_i-1}
\times\prod\{S(\x_k,\x_\l):N>k>\l,\ \s\inv(k)<\s\inv(\l)\}.\] 
In these integrals we make the substitution
\be \x_N\to{\e\ov\prod_{i<N}\x_i},\label{subs}\ee
so that $\e$ runs over a circle of radius $r^N$. The integrand becomes

\be(-1)^{|B|}\prod_{i\in B}{p+q\e\prod_{\l\ne i,\,N}\x_\l\inv-
\e\,\prod_{\l\ne N}\x_\l\inv\ov p+q\e\prod_{\l\ne i,\,N}\x_\l\inv-
\x_i}\label{first}\ee
\be\times\ \e^{x_{\s\inv(N)}-y_N-1}\,\prod_{i<N}\x_i^{x_{\s\inv(i)}-x_{\s\inv(N)}+y_N-y_i-1}\label{second}\ee
\be \times\ \prod\{S(\x_k,\x_\l):N>k>\l,\ \s\inv(k)<\s\inv(\l)\}.\label{third}\ee
The reason we still have $-1$ in the exponents in (\ref{second}) 
is that $d\x_N=\prod_{i<N}\x_i\inv\,d\e$.

\vs{1}
\noi{\bf Lemma 1}. When $|B|=1$ we have $I(\s)=0$ for all $\s\in\S_N(B)$.
\vs{1}
\noi{\bf Proof}. There is a single $i\in B$ and (\ref{first})  is analytic inside the $\x_i$-contour except for a simple pole at $\x_i=0$. The exponent of $\x_i$ in (\ref{second}) is positive since $N>i$, and so $y_N>y_i$, and $\si(i)>\si(N)$, and so $x_{\si(i)}>x_{\si(N)}$. Therefore the integrand is analytic inside the $\x_i$-contour, so the integral is zero.
   
\vs{1}
\noi{\bf Lemma 2}. When $|B|>1$ we have for all $\s\in\S_N(B)$, 
\[I(\s)=\sum_{(i,j)}I_{B,\,(i,j)}(\s),\]
where in the sum $(i,j)$ runs over all unordered pairs with $i,\,j\in B$ and $i\ne j$, where each $I_{B,\,(i,j)}(\s)$ is a lower-order integral in which (\ref{first}) is replaced by a factor depending only on $B$ and $(i,j)$ (the other factors remaining the same), and where $\x_i=\x_j$ in the domain of integration.
\vs{1}
\noi{\bf Proof}. We may assume that $q\ne0$. This case follows by a limiting argument. We are going to shrink some of the $\x_i$-contours with $i\in B$. Due to the defining property of $r$, the only poles we pass will come from the product (\ref{first}). In fact, to avoid double poles later we take $\x_i\in \C{r_i}$ with the $r_i$ all slightly different.

Take $j=\max B$ and shrink the $\x_j$-contour. The product 
(\ref{first}) has a simple pole at $\x_j=0$ (the $j$-factor has the pole and the $i$-factors with $i\ne j$ are analytic there) and the power of $\x_j$ in (\ref{second}) is positive as before, so the integrand is analytic at $\x_j=0$. For each $k\in B$ with $k\ne j$ we pass the pole at
\be\x_j={q\e\prod_{\l\ne j,k,N}\x_\l\inv\ov \x_k-p}\label{polej}\ee
coming from the $k$-factor in (\ref{first}). (Our assumption on the $r_i$ assures that there are no double poles.) For the residue we replace the $k$-factor by
\be -{p+q\e\prod_{\l\ne k,\,N}\x_\l\inv-
\e\,\prod_{\l\ne N}\x_\l\inv\ov q\e\x_j^{-2}\prod_{\l\ne j,k,N}\x_\l\inv},\label{k}\ee
where in this and the $j$-factor we replace $\x_j$ by the right side of (\ref{polej}). When $i\ne j,k$ the $i$-factor becomes
\[p+q\e\prod_{\l\ne i,\,N}\x_\l\inv-
\e\,\prod_{\l\ne N}\x_\l\inv\ov p\,(1-\x_i\x_k\inv),\]
and we replace $\x_j$ in the numerator by  the right side of (\ref{polej}).

We now shrink the $\x_k$-contour. There is a pole of order 2 at 
$\x_k=0$ coming from (\ref{k}) and the $j$-factor in (\ref{first}). Since $k<j=\max B<N$, we have $y_N-y_k\ge 2$, so the exponent of 
$\x_k$ in (\ref{second}) is at least 2. Therefore the integrand is analytic at $\x_k=0$. The factor (\ref{k}) has no other poles inside $\C{r_k}$. An $i$-factor with $i\ne j,k$ will have a pole at $\x_k=\x_i$  if $r_i<r_k$.  There is also the pole at
\[\x_k={q\e\prod_{\l\ne j,k,N}\x_\l\inv\ov \x_j-p}\]
coming from the $j$-factor. This relation and (\ref{polej}) imply $\x_j=\x_k$. 

Thus when we shrink the $\x_j$-contour and the $\x_k$-contours with $k\ne j$ we obtain $(N-2)$-dimensional integrals in each of which two of the $\x$-variables corresponding to indices in $B$ are equal. This proves the lemma. 
\vs{1}
\noi{\bf Lemma 3}. In the notation of Lemma 2, for each $(i,j)$ there is a partition of $\S_N(B)$ into pairs $\s,\,\s'$ such that $I_{B,\,(i,j)}(\s)+I_{B,\,(i,j)}(\s')=0$ for each pair.
\vs{1}
\noi{\bf Proof}. We pair $\s$ and $\s'$ if they agree except for the positions of $i$ and $j$, which are interchanged. The factor (\ref{second}) is clearly the same for both when $\x_i=\x_j$, and we shall show that the $\s$- and $\s'$-factors in (\ref{third}) are negatives of each other when $\x_i=\x_j$.

Assume for definiteness that
\be i<j\ \ {\rm and}\ \ \s\inv(i)<\s\inv(j).\label{assume}\ee
(Otherwise we reverse the roles of $\s$ and $\s'$.) Then the factor 
$S(\x_j,\x_i)$ does not appear for $\s$ in (\ref{third}) but it does appear for $\s'$.  
This factor equals $-1$ when $\x_i=\x_j$. 

To complete the proof it is enough to show that for any $k\ne i,j$ the product of $S$-factors involving $k$ and either $i$ or $j$ is the same for $\s$ and $\s'$ when $\x_i=\x_j$. If $\si(k)$ is outside the interval $(\si(i),\,\si(j))$ the $S$-factors in question are the same for $\s$ and $\s'$, so we assume $\si(k)$ is inside the interval. There are three cases, with the results displayed in the table below. The first column gives the position of $k$ relative to $i$ and $j$, the second column gives the product of $S$-factors involving $k$ and either $i$ or $j$ for $\s$, and the third column gives the corresponding product for $\s'$. 

\[\begin{array}{llllll}

i<k<j && 1 && S(\x_j,\x_k)\,S(\x_k,\x_i)\\
k<i && S(\x_i,\x_k) && S(\x_j,\x_k)\\ 
k>j && S(\x_k,\x_j) && S(\x_k,\x_i)
\end{array}\]
In all cases but the second the $S$-factors are exactly the same for $\s$ and $\s'$ when $\x_i=\x_j$. For the second we use $S(\x,\x_k)\,S(\x_k,\x)=1$. 
\vs{1}
Clearly (\ref{IBsum}) follows from Lemmas 1--3, and this completes the proof of (\ref{prob1}).

\bc {\bf III. ASEP with Multiple Species} \ec

\bc {\bf A. Bethe Ansatz solution}\ec

Observe that an interchange of particles at positions $x_i$ and $x_{i+1}$ has the same effect as leaving the particles as they were but interchanging $\p_{i}$ and $\p_{i+1}$. Thus $X$ remains the same but $\p$ is replaced by $\Tp$. (This is the same $T_i$ as before, but applied to $\p$ rather than $\s$.)

Write $u^\p(X;t)$ for $P_\Y(\X;\,t)$. The master equation for $u^\p$ differs from equation (\ref{eq1}) for $u$ since particles are not blocked as much, so there are other terms on the right side.

We define
\[\a_i(\p)=\left\{\begin{array}{ll}0&{\rm if}\ \p_{i}=\p_{i+1}\\
p&{\rm if}\ \p_{i}<\p_{i+1}\\q&{\rm if}\ \p_{i}>\p_{i+1},\end{array}\right.\]
and define $\b_i(\p)$ as above but with $p$ and $q$ interchanged. (Thus $\b_i(\p)=\a_i(\Tp)$.) 
We compute that what must be added to the right side of (\ref{eq1}) is
\[\sum_{i=1}^{N-1}\lb\a_i(\p)\,u^\Tp(x_i,x_{i+1})-\b_i(\p)\,
u^\p(x_i,x_{i+1})\rb\,\d(x_{i+1}-x_i-1).\]
(As before, we displayed entries $i$ and $i+1$.) The terms here may be incorporated into the boundary conditions, as the $\d$-terms in (\ref{eq1}) were. We conclude that if $u^\p(X;t)$ satisfies the equation
\be {\partial u^\p\ov\partial t}=\sum_{i=1}^N
\lb p\,u^\p(x_i-1)+q\,u^\p(x_i+1)-u^\p(x_i)\rb,\label{eq3}\ee
the boundary conditions for $i=1,\ld,N-1$
\[ p\,u^\p(x_i,x_i)+q\,u^\p(x_i+1,x_i+1)-u^\p(x_{i},x_{i}+1)\]
\be-\a_i(\p)\,u^\Tp(x_i,x_i+1)+\b_i(\p)\,u^\p(x_i,x_i+1)=0,\label{bc2}\ee
and the initial condition
\[u^\p(X;t)=\d_Y(X)\,\d_\n(\p),\]
then $P_\Y(\X;\,t)=u^\p(X;t)$.   

We assume a solution of the form
\[u^\p(X;t)=\sum_{\s\in \S_N} {1\ov(2\pi i)^N}\int_{\C{r}^N} A_\s^\p(\x)\,\prod_i\xs{i}^{x_i}\,\prod_i \Big(\x_i^{-y_i-1}\,e^{\ep(\x_i)\,t})\;d^N\x,\]
where $A_{e}^\p=\d_\n(\p)$ so that the initial condition is satisfied by the $\s=e$ summand. 

Substituting the left side of (\ref{bc2}) into the part of the integrand that depends on $x_i$ we now get
$(\xs{i}\,\xs{i+1})^{x_i}$ times
\[[f(b,a)+\b_i(\p)\,b]\,A_\s^\p-\a_i(\p)\,b\,A_\s^\Tp,\]
where
\be f(\x,\x')=p+q\x\x'-\x,\ \ \ a=\xs{i},\ \ b=\xs{i+1}.\label{fab}\ee
If we replace $\s$ by $\Ts$ nothing else in the integrand changes, so a sufficient condition that the integral be zero is that 
\[[f(b,a)+\b_i(\p)\,b]\,A_\s^\p+
[f(a,b)+\b_i(\p)\,a]\,A_\Ts^\p
-\a_i(\p)\,b\,A_\s^\Tp-\a_i(\p)\,a\,A_\Ts^\Tp=0.\]
If this is to hold for all $\p$ it must hold for $\p$ replaced by $\Tp$, so
\[[f(b,a)+\a_i(\p)\,b]\,A_\s^\Tp+
[f(a,b)+\a_i(\p)\,a]\,A_\Ts^\Tp
-\b_i(\p)\,b\,A_\s^\p-\b_i(\p)\,a\,A_\Ts^\p=0.\]

If $\p_{i+1}=\p_i$ then we obtain simply
\be f(b,a)\,A_\s^\p+f(a,b)\,A_\Ts^\p=0.\label{pi=pi+1}\ee
If $\p_{i+1}\ne \p_i$ we use $\b_i(\p)=1-\a_i(\p)$ to help with the computation and find by eliminating $A_{\Ts}^\Ts$ from the two equations above that
\[ A_\Ts^\p=\a_i(\p)\,{b-a\ov f(a,b)}\,A_\s^\Tp-
\lb\a_i(\p)\,{a-b\ov f(a,b)}+1\rb\,A_\s^\p.\]
(Here (\ref{fab}) was also used.) If we replace the second $\a_i(\p)$ by $1-\b_i(\p)$ we get
\[ A_\Ts^\p=\a_i(\p)\,{b-a\ov f(a,b)}\,A_\s^\Tp-
{f(b,a)-\b_i(\p)\,(a-b)\ov f(a,b)}\,A_\s^\p,\]
which now agrees with (\ref{pi=pi+1}) when $\p_{i+1}=\p_i$.

The relation is nicer for the quantity $h_\s^\p$ defined by 
\be A_\s^\p= h_\s^\p\,A_\s,\label{Ah}\ee
where $A_\s$ is defined by (\ref{A}) as before. The condition $A_{e}^\p=\d_\n(\p)$ becomes $h_{e}^\p=\d_\n(\p)$.

If we use $a-b=f(b,a)-f(a,b)$, divide by $A_\Ts=\,S(b,a)\,A_\s$, and recall (\ref{fab}), our formula becomes
\be h_\Ts^\p=h_\s^\p+\(1+S(\xs{i},\xs{i+1})\)\,\big[\a_i(\p)\,
h_\s^\Tp-\b_i(\p)\,h_\s^\p\big].\label{h}\ee
 
We have to show that these formulas together with $h_{e}^\p=\d_\n(\p)$ define $h_\s^\p$ consistently. This means that if we obtain $\s$ from $e$ by a sequence of operations $T_i$ on $\S_N$ then what we obtain from this same sequence acting on the $h_\s^\p$ is independent of the particular sequence chosen. This follows from the relations
\[{\rm(i)}\ T_i\,T_i=I;\ \ {\rm(ii)}\ T_i\,T_j=T_j\,T_i\ 
\textrm{when}\ |i-j|>1;\ \ {\rm(iii)}\ T_i\,T_{i+1}\,T_i=T_{i+1}\,T_i\,T_{i+1}.\]

We should be more precise as to what exactly these $T_i$ are. Let $\H_0$ be the set of all functions  
\[h:\p\to\textrm{function of}\ \x,\]
and $\H=\S_N\times \H_0$. Define
$T_i^0:\H\to\H_0$ by 
\be T_i^0(\s,h)=h+\(1+S(\xs{i},\xs{i+1})\)\,\big[\a_i\cdot
(h\circ T_i)-\b_i\cdot h\big].\label{T0}\ee
Then define $T_i:\H\to\H$ by
\be T_i(\s,h)=(T_i\s,T_i^0(\s,h)).\label{Ti}\ee
(The same sympol $T_i$ is used to denote these operators and operators on $\S_N$. The context should make it clear which is meant.)  It is the $T_i$ defined by (\ref{Ti}) that satisfy relations (i)--(iii).
Once we have these relations it follows from the Coxeter presentation of $\S_N$ \cite[Sec 1.9]{H} that there is a homomorphism from $\S_N$ to the group of mappings from $\H$ to itself, such that each $T_iÅ$ acting on $\S_N$ goes to the corresponding $T_i$ acting on $\H$. Because this is a homomorphism, if two products of transpositions in $\S_N$ are equal then so are the corresponding products of the $T_i$ on $\H$, which is what was claimed.

Relations (i)--(iii) are the braid relations and are the consistency equations referred to in Section I. 
Relations (i) and (ii) can be checked by hand.  Since (iii) involves only three consecutive indices, it is enough to check the case $N=3$. This was verified by a computer computation.
  
\bc{\bf B. Satisfying the initial condition}\ec

We have to show that 
\be\sum_{\s\in\S_N(B)}I^\pi(\s)=0,\label{Ipsum}\ee
where $I^\p(\s)$ is the $\s$-summand in (\ref{prob2}) with $t=0$. 

In Section II.B we made a variable change, shrank contours, and found that all $I(\s)$ with $\s\in\S_N(B)$ are sums of lower-dimensional integrals $I_{B,\,(i,j)}(\s)$, one for each unordered pair $(i,j)$ with $i,j\in B$, such that $\x_i=\x_j$ in the domain of integration. 

To extend this to the $I^\p(\s)$ we have to know first that no new poles arise from the factors $h_\s^\p$. These would come from $S$-factors involving $\x_N$, because it is only these that might introduce poles after the variable change (\ref{subs}).  Using (\ref{h}) to see which such factors could arise in the expressions for $h_\s^\p$ we start with $\s=e$ and do the following to get to our $\s\in\S_N(B)$: first, bring to the front and rearrange the $i\not\in B\cup\{N\}$. This gives some factors $1+S(\x_i,\,\x_j)$ with $i,j\ne N$. (Which exact factors occur depends on $\p$ as well as~$\s$.) Then move $N$ through the $i\in B$. Since at each step $N$ was to the right of where it moves to, we get some factors  
$1+S(\x_i,\,\x_N)$ with $i\in B$. Then rearange the elements of $B$ to reach $\s$, which introduces some factors $1+S(\x_i,\,\x_j)$ with $i,\,j\in B$. Thus the only factors involving $\x_N$ are some 
$1+S(\x_i,\,\x_N)$ with $i\in B$.\footnote
{For example, suppose $\s=(3\ 2\ 5\ 4\ 1)$, for which $B=\{4,\,1\}$. Then the steps might be
\[(1\ 2\ 3\ 4\ 5)\to(1\ 3\ 2\ 4\ 5)\to(3\ 1\ 2\ 4\ 5)\to(3\ 2\ 1\ 4\ 5)
\to(3\ 2\ 1\ 5\ 4)\to(3\ 2\ 5\ 1\ 4)\to(3\ 2\ 5\ 4\ 1).\]
The only $S$-factors involving $\x_5$ come from steps four and five, and are $S(\x_4,\x_5)$ and $S(\x_1,\x_5)$.} 
If we multiply this factor by the 
$S(\x_N,\x_i)$ from the product in (\ref{A}) for $A_\s$ we get $S(\x_N,\x_i)+1$. Thus, no new poles arise from $h_\s^\p$ with $\s\in\S_N(B)$.

We can now proceed as in Section II.B, making the variable change and shrinking contours.

We say that $\s,\,\s'$ are $(i,j)$-paired if they agree except for the positions of $i$ and $j$, which are interchanged. We saw in the proof of Lemma 3 that if $\s$ and $\s'$ are $(i,j)$-paired then the integrands in $I_{B,\,(i,j)}(\s)$ and $I_{B,\,(i,j)}(\s')$ are negatives of each other. Therefore to prove (\ref{Ipsum}) it suffices to show that if $\s$ and $\s'$ are $(i,j)$-paired then $h_\s^\p=h_{\s'}^\p$ when $\x_i=\x_j$. We show this by induction on $|\si(i)-\si(j)|$, and we may assume $\si(j)>\si(i)$.

If $\si(j)=\si(i)+1$ we replace $i$ by $\si(i)$ in (\ref{h}). Since 
$S(\x_i,\x_j)=-1$ when $\x_i=\x_j$, the statement holds then.  

Suppoes $m>1$, that the statement holds when $\si(j)=\si(i)+m-1$, and that $\si(j)=\si(i)+m$. For convenience of notation we assume that 
\[\s=(1,\, 2,\ \ld,\ m,\, m+1\ \ld),\ \ \ \s'=(m+1,\, 2,\ \ld,\ m,\ 1\ \ld).\]
(Thus $\s$ and $\s'$ are $(1,\,m+1)$-paired. Neither the actual labels nor the exact positions for the indices is relevant; only their relative positions is.) The dots represent other entries that are equal for $\s$ and $\s'$.

We want to show that $h_\s^\p=h_{\s'}^\p$ when $\x_1=\x_{m+1}$. We have 
\[T_1\s=(2,\, 1,\ \ld\ m,\, m+1,\ \ld),\ \ \ T_1\s'=(2,\, m+1,\ \ld\ m,\ 1,\ \ld).\]
Observe that $T_1\s$ and $T_1\s'$ are also $(1,\,m+1)$-paired but 
$(T_1\s)\inv(1)=2$ while $(T_1\s)\inv(m+1)=m+1$, so the induction hypothesis holds. Thus
$h_{T_1\s}^\p=h_{T_1\s'}^\p$
for all $\p$ when $\x_1=\x_{m+1}$.

Applying (\ref{h}) with $i=1$ and $\s$ replaced by $T_1\s$ and by $T_1\s'$ gives the relations

\[h_{\s}^\p=h_{T_1\s}^\p+(1+S(\x_2,\x_1))\,[\a_1(\p)\,h_{T_1\s}^{T_1\p}
-\b_1(\p)\,h_{T_1\s}^\p],\]
\[h_{\s'}^\p=h_{T_1\s'}^\p+(1+S(\x_2,\x_{m+1}))\,[\a_1(\p)\,
h_{T_1\s'}^{T_1\p}-\b_1(\p)\,h_{T_1\s'}^\p].\]
When $\x_1=\x_{m+1}$ the right sides are equal and therefore so are the left sides.

This completes the proof.

\bc{\bf C. Formulas for \boldmath$h_\s^\p$}\ec

\bc{\bf 1. A general formula}\ec

For $\s\in\S_N$ we define $h_\s\in\H_0$ to be the function $\p\to h_\s^\p$. With this notation formula (\ref{h}) may be written
$h_{T_i\s}=T_i^0(\s,h_\s)$, where $T_i^0$ is given by (\ref{T0}). From this we get 
\[h_{T_iT_j\s}=T_i^0(T_j\s,h_{T_j\s})=T_i^0(T_j\s,T_j^0(\s,h_\s))
=T_i^0 T_j(\s,h_\s).\]
And in general,
\[h_{T_{j_1}T_{j_2}\cd T_{j_m}\s}=T_{j_1}^0T_{j_2}\cd T_{j_m}(\s,h_\s).\]

If we first set $\s=e$ above and then choose the various $T_i$ such that
\be T_{j_1}T_{j_2}\cd T_{j_m}e=\s,\label{srep}\ee
we get, since $h_e^\p=\d_\n(\p)$,
\be h_\s^\p=T_{j_1}^0T_{j_2}\cd T_{j_m}(e,\d_\n).\label{hrep}\ee

For $\r\in\S_N$, and $\a$ and $\b$ functions of $\p$, we define 
\[U_i(\r,\b)=1-[1+S(\x_{\r(i)},\,\x_{\r(i+1)})]\cdot\b,\] 
\[V_i(\r,\a)=[1+S(\x_{\r(i)},\,\x_{\r(i+1)})]\cdot\a.\]
These $U_i$ and $V_i$ come from the coefficients of $h$ and $h\c T_i$, respectively, in (\ref{T0}), and products of them are the coefficients of the $\d$-summands which result from applying the various $T_j$ in (\ref{hrep}).
When we apply (\ref{T0}) consecutively  to (\ref{hrep})
we get a sum of products of the form $W_{1}\,W_{2}\cd W_m$, where each $W$ is a $U$ or $V$. 

As an example, we find that $T_j^0\,T_i(e,\d_\n)$ is a linear combination of $\d_\n,\;\d_{\n\c T_i},\;\d_{\n\c T_j}$, and $\d_{\n\c T_i\c T_j}$. The coefficient of $\d_\n$ is $U_j(T_i\,e,\b_j)\,U_i(\,e,\b_i)$, the coefficient of $\d_{\n\c T_i}$ is $U_j(T_i\,e,\b_j)\,V_i(\,e,\a_i)$, the coefficient of $\d_{\n\c T_j}$ is $V_j(T_i\,e,\a_j)\,U_i(\,e,\b_i\c T_j)$, and the coefficient of $\d_{\n\c T_i\c T_j}$ is $V_j(T_i\,e,\a_j)\,V_i(\,e,\a_i\c T_j)$. These are the four $W$-products. 

Passing to the general case, we observe that in every factor $W$ the original permutation $e$ has been composed with all earlier (rightward) $T_j$ in (\ref{hrep}). Every factor $V_i$ comes from the corresponding $h\c T_i$ summand in (\ref{T0}). This changes the $\a$ or $\b$ in each earlier factor $W$ by composing it with these $T_i$. Otherwise said, each $\a$ or $\b$ in a factor $W$ is a composition with the later $T_i$. The $W$-product is determined by these (ordered) $i$, and we denote the sequence of them by $\I=\{i_1,\ld,i_n\}$. The resulting $\d$-term is $\d_{\n\c T_{i_1}\c\cd\c T_{i_n}}$. 

In general we get a linear combination of $\d$-terms of the form $\d_{\n\c T_{i_1}\c\cd\c T_{i_n}}$, where $\I=\{i_1,\cd,i_n\}\subset \I_\s$ and $\I_\s=\{j_1,\ld,j_m\}$ is the sequence in (\ref{srep}). Each $i_k$ is a $j_\l$ while the other $j_\l$ appear between some consecutive $i_{k-1}$ and $i_k$. (For this we set $i_0=0,\;i_{n+1}=\infty$.) From the above remarks we see that the coefficient of $\d_{\n\c T_{i_1}\c\cd\c T_{i_n}}$ is
\be W_{\I}=\prod_{\l=1}^m W_{\I,\,\l}\label{Wprod}\ee
where 
\be W_{\I,\,\l}=\left\{\begin{array}{ll}V_{i_k}(T_{j_{\l-1}}\cd T_{j_m}\,e,\,
\a_{i_k}\c T_{i_{k-1}}\c\cd\c T_{i_{1}})&{\rm if}\ j_\l=i_k,\\&\\
U_{j_\l}(T_{j_{\l-1}}\cd T_{j_m}\,e,\,
\b_{j_\l}\c T_{i_{k-1}}\c\cd\c T_{i_{1}})&{\rm if}\ j_\l\in(i_{k-1},\,i_{k}).\end{array}\right.\label{W}\ee 
This gives  
\be h_\s^\p=\sum_{\I\subset\I_\s}
W_{\I}(\p)\,\d_{\n\c T_{i_1}\c\cd\c T_{i_n}}(\p),\label{hform}\ee
where $\I=\{i_1,\ld,i_n\}$.

The sequence $\I_\s$ in (\ref{srep}) is not unique, so we do not yet have an explicit formula for $h_\s^\p$. To find a particular sequence let us see how to get from $\s$ to $e$ by a sequence of $T_i$. We may do this by first bringing 1 to slot 1 by a sequence of transpositions, then bringing 2 to slot 2 by a sequence of transpositions, etc.. We first bring 1 to slot 1 by the product
\[T_1\,T_2\cd T_{\si(1)-1}.\]
Suppose we have brought $1,\ld,k-1$ to their slots. Then $k$ itself has been moved to the right by
\[\i(k):=\#\{j:j<k,\ \si(j)>\si(k)\}=\textrm{the number of inversions of the form}\ (k,j),\]
because bringing  each such $j$ to its slot has moved $k$ one slot to the right. We then bring $i$ to its slot by the product
\[T_k\cd T_{\si(k)+\i(k)-1}.\]
Thus  
\[\prod_{k=1}^{N-1}T_i\cd T_{\si(k)+\i(k)-1}\,\s=e,\]
the factor with the larger $k$ being to the left. This gives
\be\s=\prod_{k=1}^{N-1}T_{\si(k)+\i(k)-1}\cd T_k\,e,\label{srep1}\ee
the factor with the larger $k$ being to the right. With this representation we have
\be \I_\s=\bigcup_{k=1}^{N-1}\ (\l_k,\,\cd,k),\label{indices}\ee
where $\l_k=\si(k)+\i(k)-1$, the interval with the larger $k$ being to the right. 

Together with (\ref{hform}) this gives an explicit, albeit complicated, formula for $h_\s^\p$.

\bc{\bf 2. One second-class particle}\ec

Because all the $\p$ now will have 1 in a single position and 2 in the others, most of the $\a$ and $\b$ terms in (\ref{W}) will be zero when applied to $\p$. So most of the $W_{\I_\l}$ will be 0 (when a $V_j$) or 1 (when a $U_j$). It is possible systematically to determine which $W_\I$ will be nonzero when $\I_\s$ is given, and that makes it practical write down formulas for the $h_\s^\p$ in specific cases. We shall not describe the procedure but state two results that use  (\ref{indices}), which comes from representation (\ref{srep1}).

For $\n\inv(1)=1$ (the second-class particle initially in position 1) and $\p\inv(1)=j$ (the second-class particle ending in position $j$) $h_\s^\p=0$ when $\si(1)<j$ and
\[h_\s^\p=\(p-q\,S(\x_1,\xs{j})\)\,q\(1+S(\x_1,\xs{j-1})\)\cd q\(1+S(\x_1,\xs{1})\),\]
when $\si(1)\ge j$. (When $\si(1)=j$ the factor on the left equals 1.)

For $\n\inv(1)=2$ and $\p\inv(1)=j$ the formula is more complicated. When $\si(1)\ge j$ and $\si(2)+\i(2)\ge j$ it is
\[h_\s^\p=\sum_{i=1}^{j-1}\Big[\prod_{k=1}^{i-1} 
q\,\(1+S(\x_2,\xs{k})\)\cdt\(p-q\,S(\x_2,\xs{i})\)\]
\[\times \(q-p\,S(\x_1,\xs{i})\)\cdt\prod_{k=i+1}^{j-1} q\,\(1+S(\x_1,\xs{k})\)\cdt \(p-q\,S(\x_1,\xs{j})\)\Big]\]
\[+\prod_{k=1}^{j-1} q\,\(1+S(\x_2,\xs{k})\)\cdt
\(p-q\,S(\x_2,\xs{j})\)\cdt
p\,\(1+S(\x_1,\xs{j})\).\]

We should point out that in any given case (\ref{srep1}) may not be the best representation for computation. For example, take $\n\inv(1)=\p\inv(1)=3$ and $\s=(4\ 3\ 2\ 1)$. Using the representation
\be \s=T_3\,T_2\,T_1\,T_3\,T_2\,T_3\,e\label{s1}\ee
from (\ref{srep1}), there are five nonzero summands in (\ref{hform}). Therefore $h_\s^\p$ is given as a sum of five products. But if instead we use the representation 
\be\s=T_1\,T_2\,T_1\,T_3\,T_2\,T_1\,e,\label{s2}\ee
there are only two nonzero summands and we get the relatively simple formula
\[h_\s^\p=q\,(1+S(\x_2,\x_4))\,\(q-p\,S(\x_2,\x_3)\)\,q\,(1+S(\x_1,\x_3))\]
\[+\(q-p\,S(\x_2,\x_4)\)\,\(p-q\,S(\x_1,\x_4)\)\,
\(q-p\,S(\x_1,\x_3)\).\]

For another example take $\s=(4\ 3\ 2\ 1)$ as before but $\n\inv(1)=3,\,\p\inv(1)=4$. Using (\ref{s1}) we again get five nonzero summands in (\ref{hform}). But using (\ref{s2}) we get only one, and the much simpler answer
\[h_\s^\p=q\,(1+S(\x_1,\x_4))\,\(q-p\,S(\x_1,\x_3)\).\]

So although the representation (\ref{srep1}) leads to a general procedure for computing $h_\s^\p$, in any particular case there may well be a better representation for computation.

\begin{center}{\bf Acknowledgment}\end{center}

This work was supported by the National Science Foundation through grants DMS-0906387 (first author) and DMS-0854934 (second author).

\end{document}